\documentclass[reqno]{amsart} 

\pdfoutput=1
\usepackage{graphicx}
\usepackage{amssymb, amsmath, amsthm}
 \usepackage[bookmarks=false]{hyperref}

\newcommand*{\fplus}{\genfrac{}{}{0pt}{}{}{+}}
\newcommand*{\fdots}{\genfrac{}{}{0pt}{}{}{\cdots}}

\makeatletter
\renewcommand{\@biblabel}[1]{#1.} 
\makeatother

\newtheorem{Theorem}[equation]{Theorem}

\newtheorem*{rem}{Remark} 

\newtheorem*{ex}{Exercise}

\newtheorem*{rems}{Remarks} 

\newtheorem*{sol}{Solution}

\newtheorem*{nt}{Notes}

\makeatletter
\renewcommand{\@biblabel}[1]{#1.}
\makeatother

\newcommand{\qrfac}[2]{{\left({#1}; q\right)_{#2}}} 
\newcommand{\pqrfac}[3]{{\left({#1};#3\right)_{#2}}}

\numberwithin{equation}{section}

\begin{document} 
\title{How to Prove
Ramanujan's $q$-Continued Fractions }

\author{Gaurav Bhatnagar
}
\address{Educomp Solutions Ltd.}
\email{bhatnagarg@gmail.com}




\begin{abstract}
By using Euler's approach of using Euclid's algorithm to expand a power series into a continued fraction, we show how to derive Ramanujan's $q$-continued fractions in a systematic manner.
\end{abstract}


\maketitle

\begin{quote}
{\small\sl
(Ramanujan's) mastery of continued fractions was, on the formal side at any rate, beyond that of any mathematician in the world...

{\hfill $\sim$G.~H.~Hardy  \cite[p.~XXX]{hardy}}
}
\end{quote}

\section{Introduction}
The $q$-generalization of 
$$1+1+1+\cdots+1=n$$
is 
$$1+q+q^2+\cdots+q^{n-1}=\frac{1-q^n}{1-q}.$$
Similarly, Ramanujan generalized the continued fraction
\begin{equation*}\label{golden-mean-cfrac}
1+\frac{1}{1}\fplus \frac{1}{1}\fplus\frac{1}{1}\fplus\fdots
\end{equation*}
to
\begin{equation*}
1+\frac{q}{1}\fplus\frac{q^2}{1}\fplus\frac{q^3}{1}\fplus\fdots,
\end{equation*}
and showed that for $|q|<1$, this continued fraction can be written as a ratio of very similar looking sums. 
\begin{align}
{1}+ \frac{q}{1}\fplus \frac{q^2}{1}\fplus\frac{q^3}{1}\fplus\fdots &=
\frac{\displaystyle\sum_{k=0}^{\infty} \frac{q^{k^2}}{(1-q)(1-q^2)\cdots (1-q^k)}}
{\displaystyle\sum_{k=0}^{\infty} \frac{q^{k^2+k}}{{(1-q)(1-q^2)\cdots (1-q^k)}}} .\label{rr-cfrac1}
\end{align}


From such a humble beginning, Ramanujan wrote down several generalizations and special cases, in the process sometimes rediscovering some continued fractions found earlier by Gauss, Eisenstein and Rogers.  As was his way, he did not record his proofs. 

Proofs were provided over the years, by many mathematicians.  We mention specially Andrews \cite{andrews-ln1} and Adiga, Berndt,  Bhargava, and Watson \cite{abbw}.  Proofs have been compiled in
\cite{LN1} and 
 \cite{berndt-notebooks-3}.  

The purpose of this tutorial  is to help the reader learn how to prove some of Ramanujan's  $q$-continued fraction formulas. 
In particular,  we will show how to derive 
nine continued fractions that appear in
Ramanujan's Lost Notebook \cite{LN} (see 
Andrews and Berndt \cite[ch.~6]{LN1}) and his earlier Second Notebook (see Berndt \cite[ch.~16]{berndt-notebooks-3}).  Lorentzen~\cite{ll98} has provided an alternative approach to Ramanujan's earlier continued fractions not covered here. 

Ramanujan was a master of manipulatorics in the class of Euler himself. Thus  it is  appropriate that  the continued fraction formulas of Ramanujan here are all derived by using the same approach as the one taken by  Euler \cite{eulerE616}  
for the \lq\lq {\em transformation of the divergent series 
$1-mx+m(m+n)x^2
-m(m+n)(m+2n)x^3+m(m+n)(m+2n)(m+3n)x^4+$ etc. 
into a continued fraction}".

\section{Euler's approach}\label{sec:euler}

Euler~\cite{eulerE616} used the elementary identity:
\begin{equation}\label{div-1step}
\frac{N}{D}=1+\frac{N-D}{D}.
\end{equation}
This is  simply one step of long division, provided the quotient when the numerator $N$ is divided by the denominator $D$ is $1$. This identity is used to \lq divide\rq\ a formal power series of the form $1+a_1z+a_2z^2+\cdots$  with another series of the same form. This elementary identity is used repeatedly to derive all the continued fraction formulas in this paper.


The following exercise will help you discover for yourself the key idea of Euler's approach. 
\begin{ex}
Use \eqref{div-1step} repeatedly to expand the fraction $13/8$ into a continued fraction.
\end{ex}
You may also enjoy spotting \eqref{div-1step} (and a continued fraction mentioned in the introduction) in Gowers' \cite[p.~41--45]{gowers} proof that  the Golden Ratio is irrational.  

\section{The Rogers-Rama\-nujan Continued Fraction}\label{sec:RR}

We proceed to apply Euler's approach to prove a slight generalization of \eqref{rr-cfrac1}, the famous Rogers-Ramanujan continued fraction. This continued fraction is due to Rogers \cite{rogers} and Ramanujan  \cite[ch.~16]{berndt-notebooks-3}. 

All the continued fractions considered in this paper have a special parameter $q$ in them. The associated series found here are of a particular type known as $q$-series. They are recognizable by the presence of the  $q$-rising factorial $\qrfac{a}{k}$, defined as:
$$\qrfac{a}{k} :=
\begin{cases}
1 &{\text{ if } k=0},\\
(1-a)(1-aq)\cdots (1-aq^{k-1}) &{\text{ if }} k\geq 1.\\
\end{cases}
$$
Similarly, the
infinite $q$-rising factorial is defined as:
$$\qrfac{A}{\infty}:=\prod_{j=0}^{\infty} (1-Aq^j), {\text{ for $|q|<1$}}.$$

The Rogers-Ramanujan continued fraction \cite[Cor.\ to Entry 15, ch.\ 16]{berndt-notebooks-3} is \eqref{rr-cfrac1} with one additional parameter:
\begin{equation}\label{cor-entry15}
\frac{\displaystyle\sum_{k=0}^{\infty} \frac{q^{k^2+k}}{\qrfac{q}{k}}a^k}
{\displaystyle\sum_{k=0}^{\infty} \frac{q^{k^2}}{\qrfac{q}{k}}a^k}
=
\frac{1}{1}\fplus \frac{aq}{1}\fplus\frac{aq^2}{1}\fplus\fdots.
\end{equation}
Its derivation is as follows. 

The first step is to rewrite the ratio of sums on the LHS of \eqref{cor-entry15} as 
\begin{align*}
\frac{\displaystyle\sum_{k=0}^{\infty} \frac{q^{k^2+k}}{\qrfac{q}{k}}a^k}
{\displaystyle\sum_{k=0}^{\infty} \frac{q^{k^2}}{\qrfac{q}{k}}a^k}
&= \frac{1}{\frac{\displaystyle\sum_{k=0}^{\infty} \frac{q^{k^2}}{\qrfac{q}{k}}a^k}
{\displaystyle\sum_{k=0}^{\infty} \frac{q^{k^2+k}}{\qrfac{q}{k}}a^k}}.
\end{align*}
Next, formally \lq divide\rq\    the two sums in the denominator by using the elementary identity \eqref{div-1step},
to obtain:  
\begin{equation*}
\frac{1}{1+ \frac{\displaystyle\sum_{k=0}^{\infty} \frac{q^{k^2}}{\qrfac{q}{k}}a^k-
\displaystyle\sum_{k=0}^{\infty} \frac{q^{k^2+k}}{\qrfac{q}{k}}a^k}
{\displaystyle\sum_{k=0}^{\infty} \frac{q^{k^2+k}}{\qrfac{q}{k}}a^k}}.
\end{equation*}
Now consider the difference of sums
$$\displaystyle\sum_{k=0}^{\infty} \frac{q^{k^2}}{\qrfac{q}{k}}a^k-
\displaystyle\sum_{k=0}^{\infty} \frac{q^{k^2+k}}{\qrfac{q}{k}}a^k
=\sum_{k=0}^{\infty} \frac{q^{k^2}}{\qrfac{q}{k}}a^k(1-q^k),$$
where we have subtracted the two sums term-by-term. Note that the first term (corresponding to the index $k=0$) is $0$, and the sum actually runs from $k=1$ to $\infty$. But
\begin{align*}
\sum_{k=1}^{\infty} \frac{q^{k^2}}{\qrfac{q}{k}}a^k(1-q^k) 
&=
\sum_{k=1}^{\infty} \frac{q^{k^2}}{\qrfac{q}{k-1}}a^k \cr
&=\sum_{k=0}^{\infty} \frac{q^{(k+1)^2}}{\qrfac{q}{k}}a^{k+1} \cr
&=aq\sum_{k=0}^{\infty} \frac{q^{k^2+2k}}{\qrfac{q}{k}}a^{k}.
\end{align*}
We have shifted the index so that the sum once again runs from $0$ to $\infty$. In the process, we extracted the common factor $aq$ from the sum. 
The ratio of sums on the LHS of \eqref{cor-entry15} can be written as
\begin{equation*}
\frac{1}{1+ 
{ 
\frac{\displaystyle aq \sum_{k=0}^{\infty} \frac{q^{k^2+2k}}{\qrfac{q}{k}}a^{k}}
{\displaystyle\sum_{k=0}^{\infty} \frac{q^{k^2+k}}{\qrfac{q}{k}}a^k}
}}
=\frac{1}{1} \fplus \frac{aq}
{
\frac{\displaystyle\sum_{k=0}^{\infty} \frac{q^{k^2+k}}{\qrfac{q}{k}}a^k}
{\displaystyle \sum_{k=0}^{\infty} \frac{q^{k^2+2k}}{\qrfac{q}{k}}a^{k}}
}.
\end{equation*}
Once again, divide the two sums using \eqref{div-1step} and find that the LHS of 
\eqref{cor-entry15} equals:
\begin{equation*}
\frac{1}{1} \fplus \frac{aq}
{
1+\frac{\displaystyle\sum_{k=0}^{\infty} \frac{q^{k^2+k}}{\qrfac{q}{k}}a^k(1-q^k)}
{\displaystyle \sum_{k=0}^{\infty} \frac{q^{k^2+2k}}{\qrfac{q}{k}}a^{k}}
}.
\end{equation*}
Now, as before, note that:
\begin{align*}
\sum_{k=0}^{\infty} \frac{q^{k^2+k}}{\qrfac{q}{k}}a^k(1-q^k) &=
\sum_{k=1}^{\infty} \frac{q^{k^2+k}}{\qrfac{q}{k-1}}a^k \cr
&=\sum_{k=0}^{\infty} \frac{q^{(k+1)^2+k+1}}{\qrfac{q}{k}}a^{k+1} \cr
&=aq^2\sum_{k=0}^{\infty} \frac{q^{k^2+3k}}{\qrfac{q}{k}}a^{k}.
\end{align*}
This time the common factor $aq^2$ pops out, and we find that the LHS of 
\eqref{cor-entry15} can be written as:
\begin{equation*}
\frac{1}{1} \fplus \frac{aq}{1}\fplus \frac{aq^2}
{
\frac{\displaystyle\sum_{k=0}^{\infty} \frac{q^{k^2+2k}}{\qrfac{q}{k}}a^k}
{\displaystyle \sum_{k=0}^{\infty} \frac{q^{k^2+3k}}{\qrfac{q}{k}}a^{k}}
}.
\end{equation*}
This process can be repeated. 

The pattern is clear. Define $R(s)$, for $s=0,1,2, \dots$, as follows:
\begin{equation}\label{RR-Rs}
R(s):=\sum_{k=0}^{\infty} \frac{q^{k^2+sk}}{\qrfac{q}{k}}a^k.
\end{equation}
Then, using \eqref{div-1step} we have
{\allowdisplaybreaks
\begin{align*}
\frac{R(s)}{R(s+1)}&=1+\frac{R(s)-R(s+1)}{R(s+1)}\\
&=1+ \frac{1}{R(s+1)}
{\displaystyle\sum_{k=0}^{\infty} \frac{q^{k^2+sk}}{\qrfac{q}{k}}a^k(1-q^k)}
\\
&=1+\frac{1}{R(s+1)}
{\sum_{k=1}^{\infty} \frac{q^{k^2+sk}}{\qrfac{q}{k-1}}a^k}\\
&=1+\frac{1}{R(s+1)}
{\sum_{k=0}^{\infty} \frac{q^{k^2+2k+1+sk+s}}{\qrfac{q}{k}}a^{k+1}}\\
&=1+\frac{aq^{s+1}}{R(s+1)}
{\sum_{k=0}^{\infty} \frac{q^{k^2+(s+2)k}}{\qrfac{q}{k}}a^{k}}\\
&=1+\frac{aq^{s+1}}{{\cfrac{R(s+1)}{R(s+2)}}}.
\end{align*}}%
Thus we obtain:
\begin{equation}\label{RR-recursion}
\frac{R(s)}{R(s+1)} =
1+\frac{aq^{s+1}}{{\cfrac{R(s+1)}{R(s+2)}}}.
\end{equation}
This gives, on iteration, 
\begin{equation}\label{RR-iteration}
\frac{R(1)}{R(0)}=\frac{1}{\cfrac{R(0)}{R(1)}}=
\frac{1}{1}\fplus \frac{aq}{1}\fplus \frac{aq^2}{1}\fplus \fdots\fplus \frac{aq^{s+1}}{{\cfrac{R(s+1)}{R(s+2)}}}.
\end{equation}
Take limits as $s\to \infty$ to formally obtain \eqref{cor-entry15}. We will show how to justify the limiting process in the next section.  \qed

To summarize, Euler's approach is as follows. Divide the two series using \eqref{div-1step}. Now cancel the first term, shift the index, and extract the common terms so that the difference series once again has $1$ as the first term. A few steps of this, and one can find the first few terms of the continued fraction. Its now easy to find a pattern and to prove it. 

A large percentage of continued fractions involving series (whether $q$-series or not) can be found just like this. 

The following exercise will help you come to terms with the notation of $q$-rising factorials. 
\begin{ex}
Show that, for $|q|<1$: 
\end{ex}
\begin{enumerate}
\item
$\displaystyle \qrfac{q}{\infty}\qrfac{-q}{\infty}=\pqrfac{q^2}{\infty}{q^2}.$
\item\label{ex3.2}
$\displaystyle \lim_{a\to 0}\qrfac{-\lambda/a}{k}a^k=\lambda^k q^{k\choose 2}$,
 where $k=0, 1, 2, \dots$.
\item
$\displaystyle{\pqrfac{-q^2}{\infty}{q^2}}=\frac{1}{\pqrfac{q^2}{\infty}{q^4}}.$
\item
$\displaystyle{\pqrfac{-q}{\infty}{q^2}}=\frac{\pqrfac{q^2}{\infty}{q^4}}{\pqrfac{q}{\infty}{q^2}}.$
\end{enumerate}

\section{Convergence matters}\label{sec:convergence}

Before proceeding with the derivation of more continued fractions, we record our attitude towards matters like convergence.  

While in the flow of performing symbolic calculations, one doesn't usually think about such things.
Most of the time we regard series as formal power series and  focus on the algebraic calculations.  
It is only after getting a nice formula that we worry about whether the steps can be made rigorous or not, and whether the symbolic calculations actually yield a valid formula. 

The objective of this section is to provide the necessary background information required to ensure that the symbolic calculations are indeed rigorous. We also provide links to references with more information. 

Infinite products such as $\qrfac{A}{\infty}$ converge when $|q|<1$. For an introduction to convergence of infinite products, try Rainville~\cite{rainville}. The test for convergence of infinite products says that 
$$\prod_{j=0}^{\infty} (1-Aq^j)$$ converges absolutely provided the sum
$$A\sum_{j=0}^{\infty} q^j$$ converges absolutely.  Which it does, whenever $|q|<1$. 

In the case of series, the convergence conditions follow from the ratio test. Most of the series  considered here are of the form $$\sum_{k=0}^{\infty} a_k z^k,$$ where $a_k$ has factors such as $\qrfac{a}{k}$. When we take the ratio $a_{k+1}/a_k$, only $(1-aq^k)$ remains of this. Now when $|q|<1$, this does not contribute anything to the ratio test. So many of these series will converge as long as the power series variable satisfies $|z|<1$. Occasionally, as in \eqref{cor-entry15}, we may get lucky and get a quadratic power $q^{k^2}$ in the summand of the series. In this case we don't even need any conditions on the power series variable $z$. The series converges whenever $|q|<1$. 

You may refer to Gasper and Rahman~\cite{grhyp}  for more details on the convergence conditions for the kinds of series and products considered here. In what follows, we simply state the convergence conditions without comment. 

Now we come to the convergence of continued fractions. Our exposition is based on Lorentzen and Waadeland \cite{lw92}. 

Let $\{ a_n\}_{n=1}^{\infty}$ and $\{ b_n\}_{n=0}^{\infty}$   be sequences of complex numbers, $a_n\neq 0$. Consider the sequence of mappings 
$$s_0(w)=b_0+w; s_n(w)=\frac{a_n}{b_n+w}, \text{ for } n=1, 2, \dots $$
These mappings are from $\hat{\mathbb{C}}$ to $\hat{\mathbb{C}}$, where $\hat{\mathbb{C}}$ are the extended complex numbers ${\mathbb{C}}\cup \left\{\infty\right\}$.
Let $S_n(w)$ be defined as follows:
$$S_0(w)=s_0(w)=b_0+w; S_n(w)=S_{n-1}(s_n(w)).$$
Then $S_n(w)$ can be written as
$$S_n(w)=b_0+\frac{a_1}{b_1}\fplus\frac{a_2}{b_2}\fplus\fdots\fplus\frac{a_n}{b_n+w}.$$
Note that since $a_n\neq 0$, the $s_k(w)$ are non-singular linear fractional transformations.  The $S_n$ are compositions of these, and are thus also non-singular linear fractional transformations. 
 
A continued fraction is an ordered pair
$\left( (\{ a_n\}, \{ b_n\}), S_n(0)\right)$, where $a_n$, $b_n$ and $S_n$ are as above. When $S_n(0)$ converges to an extended complex number $S$, we say that the continued fraction converges, and we write
$$S=b_0+\frac{a_1}{b_1}\fplus\frac{a_2}{b_2}\fplus\fdots.$$

We call $S_n(0)$ the {\em approximants} of the continued fraction. The $S_n(w)$ are called {\em modified approximants}. The convergence of 
$S_n(w)$ (or even $S_n(w_n)$, for suitably chosen $w_n$) is called {\em modified convergence}. Modified convergence is distinguished from {\em ordinary convergence}, the convergence of $S_n(0)$ discussed above. 

Modified convergence is an immediate consequence of Euler's approach for all the continued fractions presented in this paper.  However, we can obtain the ordinary convergence of the continued fractions too,  by appealing to a convergence theorem due to Worpitzsky. 

We require the shorthand notation 
$$\mathop{\mathbf{K}} \left( {a_n}/{b_n}\right)
\equiv
\mathop{\mathbf{K}}\limits_{n=1}^{\infty} \left( {a_n}/{b_n}\right)=
\frac{a_1}{b_1}\fplus\frac{a_2}{b_2}\fplus\fdots.$$
\begin{Theorem}[Worpitzky's Theorem {\cite[p.~35]{lw92}}] Let $| a_n|\leq 1/4$ for all $n\geq 1$. Then 
$\mathop{\mathbf{K}} \left( {a_n}/{b_n}\right)$ converges. All approximants are in the disk $|w|<1/2$, and the value of the continued fraction is in the disk  $|w|\leq 1/2.$ Moreover, the convergence of $S_n(w)$ is uniform with respect to $w$, for $|w|<1/2$. 
\end{Theorem}
Suppose that the approximants converge to $S$, and suppose that $|w_n|<1/2$. Suppose further that $w_n\to 0$. Then, since the convergence is uniform, the modified approximants $S_n(w_n)$ will converge to $S$ too. We will use Worpitzky's theorem in precisely these conditions. 

Let us now apply these ideas to discuss convergence in the proof of the Rogers-Ramanujan continued fraction \eqref{cor-entry15}. 

First we show that modified convergence follows immediately from \eqref{RR-iteration}. In the case of the Rogers-Ramanujan continued fraction, the modified approximants are as follows:
$S_0(w)=w$, and 
\begin{equation*}
S_n(w)=
\frac{1}{1}\fplus \frac{aq}{1}\fplus \frac{aq^2}{1}\fplus \fdots\fplus \frac{aq^{n-1}}{1+w}.
\end{equation*}
Further, define $w_n$ as:
$w_0:=1$, and 
$$w_n:=aq^n\frac{R(n+1)}{R(n)},$$
where $R(s)$ is given by \eqref{RR-Rs}. 
We have shown, in \eqref{RR-iteration} that for $n\geq 1$
$$S_n(w_n)=\frac{R(1)}{R(0)}.$$
Clearly, 
$$\lim_{n\to\infty} S_n(w_n)=\frac{R(1)}{R(0)}$$
and we have modified convergence. 

To show the convergence of $S_n(0)$, we need to work with the \lq tail' of the continued fraction, see \cite[Th.~1, p.~56]{lw92}. The $N$th tail of the continued fraction 
$$b_0+\mathop{\mathbf{K}} \left( {a_n}/{b_n}\right)$$  is  $$\mathop{\mathbf{K}} \left( {a_{N+n}}/{b_{N+n}}\right).$$
The continued fraction
$$b_0+\mathop{\mathbf{K}} \left( {a_n}/{b_n}\right)$$
converges if and only if its $N$th tail converges for some natural number $N$. Further, if the $n$th approximant $S_n^{(N)}(0)$ of the $N$th tail converges to $S^{(N)}$, then $S_n(0)$ converges to 
$$S=b_0+\frac{a_1}{b_1}\fplus\frac{a_2}{b_2}\fplus\fdots\fplus\frac{a_N}{b_N+S^{(N)}}.$$

Returning to the Rogers-Ramanujan continued fraction, note that if $|q|<1$, we can find an  $N$ such that 
$$\left| aq^n\right| <\frac{1}{4}\text{ for all } n\geq N$$
and
$$\left| aq^n\frac{R(n+1)}{R(n)}\right| <\frac{1}{2}\text{ for all } n\geq N.$$
We now consider the $N$th tail of the Rogers-Ramanujan continued fraction. Worpitzky's Theorem applies to this continued fraction, due to the conditions above. The $n$th approximant is given by
\begin{equation*}
S_n^{(N)}(w)=
\frac{aq^N}{1}\fplus \frac{aq^{N+1}}{1}\fplus \fdots\fplus \frac{aq^{N+n-1}}{1+w}.
\end{equation*}
Here $w_n$ are defined as:
$$w_n:=aq^{N+n}\frac{R(N+n+1)}{R(N+n)}.$$ 
By iterating \eqref{RR-recursion}, we can see that
$$S_n^{(N)}(w_n)=aq^N\frac{R(N+1)}{R(N)}.$$
By our comments following the statement of Worpitzky's Theorem, we find that this implies that the modified approximants converge to the same value as the approximants $S_n^{(N)}(0)$. 
This implies that $S_n(0)$, the $n$th approximant of the Rogers-Ramanujan continued fraction, converges to
$$\frac{1}{1}\fplus \frac{aq}{1}\fplus \frac{aq^2}{1}\fplus \fdots\fplus \frac{aq^{N}}{{\cfrac{R(N)}{R(N+1)}}},$$
which is equal to $R(1)/R(0)$ by \eqref{RR-iteration}. This completes the justification of the limiting process in the proof of \eqref{cor-entry15}.

The above arguments apply to all the continued fractions in this paper. In each case, we will find that modified convergence is an immediate consequence of our approach. To show that modified convergence implies convergence, we use Worpitzky's Theorem on the tail of the continued fraction. To do so, we will require $|q|<1$ and occasionally some further restrictions on other parameters.

To be able to apply Worpitzsky's theorem, we will often need to write one of Ramanujan's continued fractions in an equivalent form. Two continued fractions are called {\em equivalent} if they have the same sequence of approximants. 
\begin{ex} Show that the continued fraction
$$
\frac{1}{1}\fplus\frac{\lambda q}{1+bq}\fplus\frac{\lambda q^2}{1+bq^2}\fplus\frac{\lambda q^3}{1+bq^3}\fplus
\frac{\lambda q^4}{1+bq^4}\fplus\fdots 
$$
is equivalent to 
\begin{align*}
\frac{1}{1}\fplus\frac{\lambda q/{(1+bq)}}{1}\fplus&\frac{\lambda q^2/{(1+bq)(1+bq^2)}}{1}\fplus \\
&\frac{\lambda q^3/{(1+bq^2)(1+bq^3)}}{1}\fplus\frac{\lambda q^4/{(1+bq^3)(1+bq^4)}}{1}\fplus\fdots .
\end{align*}
Given that $|q|<1$, show that these continued fractions converge.
\end{ex}
This continued fraction (again due to Ramanujan) is slightly more general that the Rogers-Ramanujan continued fraction. It also appears as a ratio of two similar looking sums. 

\section{More continued fractions by Euler's approach}\label{sec:euler2}
Next, we consider, as did Ramanujan, the slightly more general
expression
\begin{align}
g(b,\lambda)&:=\sum_{k=0}^{\infty}\frac{ q^{k^2}}{\qrfac{q}{k}\qrfac{-bq}{k}}\lambda^k \label{g}\\
&=1+\frac{\lambda q}{(1-q)(1+bq)}+\frac{\lambda^2 q^4}{(1-q)(1-q^2)(1+bq)(1+bq^2)}+\cdots
.\notag
\end{align}
Ramanujan found the continued fraction (see Entry 15 of \cite[ch.\ 16 ]{berndt-notebooks-3} or \cite[Entry 6.3.1(ii)]{LN1})
\begin{align}
\frac{g(b,\lambda q)}{g(b,\lambda)}
&=
 \frac{\displaystyle
\sum_{k=0}^{\infty}\frac{q^{k^2+k}}{\qrfac{q}{k}\qrfac{-bq}{k}}\lambda^k}
{\displaystyle\sum_{k=0}^{\infty}\frac{q^{k^2}}{\qrfac{q}{k}\qrfac{-bq}{k}} \lambda^k}
\label{g-frac-sums}
\\
&=\frac{1}{1}\fplus\frac{\lambda q}{1+bq}\fplus\frac{\lambda q^2}{1+bq^2}\fplus\frac{\lambda q^3}{1+bq^3}\fplus\fdots .\label{g-cfrac2}
\end{align}
Note that when $b=0$ and $\lambda =a$, this reduces to \eqref{cor-entry15}.

%
The proof follows the pattern of the proof of the generalized Rogers-Ramanujan continued fraction presented in 
\S\ref{sec:RR}. The first step is to note that
$$\frac{g(b,\lambda q)}{g(b,\lambda)}
=
\frac{1}{
\frac{\displaystyle\sum_{k=0}^{\infty} \frac{q^{k^2}}{\qrfac{q}{k}\qrfac{-bq}{k}}\lambda^k}
{\displaystyle\sum_{k=0}^{\infty} \frac{q^{k^2+k}}{\qrfac{q}{k}\qrfac{-bq}{k}}\lambda^k}
}.
$$
Once again, we divide the two sums in the denominator using \eqref{div-1step} and find that
 \eqref{g-frac-sums} equals:
\begin{equation*}
{\frac{1}
{1}\fplus\frac{\displaystyle\sum_{k=0}^{\infty} \frac{q^{k^2}}{\qrfac{q}{k}\qrfac{-bq}{k}}\lambda^k(1-q^k)}
{\displaystyle \sum_{k=0}^{\infty} \frac{q^{k^2+k}}{\qrfac{q}{k}\qrfac{-bq}{k}}\lambda^{k}}
}.
\end{equation*}
Now, note that:
\begin{align*}
\sum_{k=0}^{\infty} \frac{q^{k^2}}{\qrfac{q}{k}\qrfac{-bq}{k}}\lambda^k(1-q^k) &=
\sum_{k=1}^{\infty} \frac{q^{k^2}}{\qrfac{q}{k-1}\qrfac{-bq}{k}}\lambda^k \cr
&=\sum_{k=0}^{\infty} \frac{q^{(k+1)^2}}{\qrfac{q}{k}\qrfac{-bq}{k+1}}\lambda^{k+1} \cr
&=\frac{\lambda q}{1+bq} \sum_{k=0}^{\infty} \frac{q^{k^2+2k}}{\qrfac{q}{k}\qrfac{-bq^2}{k}}\lambda^{k}.
\end{align*}
The common factor $\lambda q/(1+bq)$ pops out, and we find that 
\eqref{g-frac-sums} can be written as:
\begin{equation}\label{biggerCF-c1}
\frac{1}{1}\fplus \frac{\lambda q}{
(1+bq)\cfrac{
{\displaystyle \sum_{k=0}^{\infty} \frac{q^{k^2+k}}{\qrfac{q}{k}\qrfac{-bq}{k}}\lambda^{k}}
}
{\displaystyle \sum_{k=0}^{\infty} \frac{q^{k^2+2k}}{\qrfac{q}{k}\qrfac{-bq^2}{k}}\lambda^{k}
}}.
\end{equation}
Next, use \eqref{div-1step} again to find that this expression can be written as:
\begin{equation*}
\frac{1}{1}\fplus\frac{\lambda q}{
{(1+bq)\left( 1+
\frac{{\displaystyle \sum_{k=0}^{\infty} \frac{q^{k^2+k}}{\qrfac{q}{k}\qrfac{-bq}{k}}\lambda^{k}}
-\displaystyle \sum_{k=0}^{\infty} \frac{q^{k^2+2k}}{\qrfac{q}{k}\qrfac{-bq^2}{k}}\lambda^{k}}
{\displaystyle \sum_{k=0}^{\infty} \frac{q^{k^2+2k}}{\qrfac{q}{k}\qrfac{-bq^2}{k}}\lambda^{k}
}\right)}.
}
\end{equation*}
Consider the difference of sums appearing in this expression. 
\begin{align}\label{calc1-biggerCF1}
{\displaystyle \sum_{k=0}^{\infty} \frac{q^{k^2+k}}{\qrfac{q}{k}\qrfac{-bq}{k}}\lambda^{k}}
&-\displaystyle \sum_{k=0}^{\infty} \frac{q^{k^2+2k}}{\qrfac{q}{k}\qrfac{-bq^2}{k}}\lambda^{k}\cr
&=\displaystyle \sum_{k=0}^{\infty} \frac{q^{k^2+k}}{\qrfac{q}{k}\qrfac{-bq}{k}}\lambda^{k}
\left[ 1-\frac{1+bq}{1+bq^{k+1}}q^k\right]\cr
&=\displaystyle \sum_{k=0}^{\infty} \frac{q^{k^2+k}}{\qrfac{q}{k}\qrfac{-bq}{k}}\lambda^{k}
\left[ \frac{1-q^k}{1+bq^{k+1}}\right].
\end{align}
Note that the term corresponding to the index $k=0$ in \eqref{calc1-biggerCF1} is $0$, so the sum actually runs from $k=1$ to $\infty$. So we find that that this sum equals
\begin{align}\label{calc1-biggerCF2}
&\displaystyle \sum_{k=1}^{\infty} \frac{q^{k^2+k}}{\qrfac{q}{k}\qrfac{-bq}{k}}\lambda^{k}
\left[ \frac{1-q^k}{1+bq^{k+1}}\right]\cr
&=\displaystyle \sum_{k=1}^{\infty} \frac{q^{k^2+k}}{\qrfac{q}{k-1}\qrfac{-bq}{k+1}}\lambda^{k}
\cr
&=\displaystyle \sum_{k=0}^{\infty} \frac{q^{k^2+3k+2}}{\qrfac{q}{k}\qrfac{-bq}{k+2}}\lambda^{k+1}
\cr
&=\frac{\lambda q^2}{(1+bq)(1+bq^2)}
\displaystyle \sum_{k=0}^{\infty} \frac{q^{k^2+3k}}{\qrfac{q}{k}\qrfac{-bq^3}{k}}\lambda^{k}.
\end{align}
Using equation \eqref{calc1-biggerCF2}, the ratio \eqref{g-frac-sums} can be written as:
\begin{equation*}
\frac{1}{1}\fplus\frac{\lambda q}{1+bq}\fplus \frac{\lambda q^2}
{
(1+bq^2)\cfrac{
{\displaystyle \sum_{k=0}^{\infty} \frac{q^{k^2+2k}}{\qrfac{q}{k}\qrfac{-bq^2}{k}}\lambda^{k}}
}
{\displaystyle \sum_{k=0}^{\infty} \frac{q^{k^2+3k}}{\qrfac{q}{k}\qrfac{-bq^3}{k}}\lambda^{k}
}}.
\end{equation*}
The pattern is clear. It is apparent that this process will yield the continued fraction \eqref{g-cfrac2}. 

We define $g_1(s)$, for $s=0,1, 2, 3, \dots$, as follows:
\begin{equation}\label{g2s}
g_1(s):={\displaystyle\sum_{k=0}^{\infty} \frac{q^{k^2+sk}}{\qrfac{q}{k}\qrfac{-bq^{s}}{k}}\lambda^k}.
\end{equation}

Now, using \eqref{div-1step} we have, for $s=0,1, 2, 3,\dots$:
\begin{align*}\allowdisplaybreaks
\frac{g_1(s)}{g_1(s+1)}&=1+\frac{g_1(s)-g_1(s+1)}{g_1(s+1)}\cr
&=1+ \frac{1}{g_1(s+1)}
{\displaystyle\sum_{k=0}^{\infty} \frac{q^{k^2+sk}}{\qrfac{q}{k}\qrfac{-bq^{s}}{k}}\lambda^k
\left[1-\frac{1+bq^s}{1+bq^{s+k}}q^k\right]}
\cr
&=1+\frac{\lambda q^{s+1}}{(1+bq^s)(1+bq^{s+1})g_1(s+1)}
{\sum_{k=0}^{\infty} \frac{q^{k^2+(s+2)k}}{\qrfac{q}{k}\qrfac{-bq^{s+2}}{k}}\lambda^{k}}\cr
&=1+\frac{\lambda q^{s+1}}{(1+bq^s)(1+bq^{s+1}){\cfrac{g_1(s+1)}{g_1(s+2)}}},
\end{align*}
where we have skipped a few steps similar to \eqref{calc1-biggerCF1} and \eqref{calc1-biggerCF2} above. 

Multiplying both sides by $(1+bq^s)$, we obtain, for $s=0,1, 2, 3, \dots$:
\begin{equation}\label{g2-recursion}
(1+bq^s)\frac{g_1(s)}{g_1(s+1)}=
1+bq^s+\frac{\lambda q^{s+1}}{(1+bq^{s+1}){\cfrac{g_1(s+1)}{g_1(s+2)}}}.
\end{equation}

We have already shown in \eqref{biggerCF-c1} that \eqref{g-frac-sums} is equal to
$$
\frac{g(b,\lambda q)}{g(b,\lambda)}=
\frac{1}{1}\fplus\frac{\lambda q}{(1+bq)\cfrac{g_1(1)}{g_1(2)}}.$$
Now by using the $s=1, 2, 3, \dots$ case of \eqref{g2-recursion} 
we find that: 
\begin{equation}\label{biggerCF1-c2}
\frac{g(b,\lambda q)}{g(b,\lambda)}=
\frac{1}{1}\fplus\frac{\lambda q}{1+bq}\fplus \frac{\lambda q^2}{1+bq^2}\fplus\fdots 
\fplus\frac{\lambda q^{s+1}}{(1+bq^{s+1})\cfrac{g_1(s+1)}{g_1(s+2)}}.
\end{equation}
Now by arguments of \S \ref{sec:convergence}, we can show modified convergence of the continued fraction, as $s\to \infty$.  To show ordinary convergence, we apply Worpitzky's Theorem to an equivalent continued fraction, given in the Exercise of \S \ref{sec:convergence}. In this manner,  we complete the proof of \eqref{g-cfrac2}.

Its time now to try your hand at using Euler's approach to derive a continued fraction. 
 Consider the sum:
\begin{equation}\label{G}
G(a,b,\lambda):=\sum_{k=0}^{\infty}
\frac{\qrfac{-\lambda/a}{k} q^{(k^2+k)/2}}{\qrfac{q}{k}\qrfac{-bq}{k}}a^k.
\end{equation}
\begin{ex}
Show that, for $|q|<1$:
\begin{align}
&\frac{G(aq, b, \lambda q)}{G(a, b, \lambda)}
=
 \frac{\displaystyle
\sum_{k=0}^{\infty}\frac{\qrfac{-\lambda/a}{k}q^{(k^2+3k)/2}}{\qrfac{q}{k}\qrfac{-bq}{k}}a^k}
{\displaystyle\sum_{k=0}^{\infty}\frac{\qrfac{-\lambda/a}{k}q^{(k^2+k)/2}}{\qrfac{q}{k}\qrfac{-bq}{k}}a^k}
\label{G-frac-sums}
\\
&=
\frac{1}{1}\fplus\frac{aq+\lambda q}{1+bq}\fplus\frac{\lambda q^2-abq^3}
{1+bq^2}\fplus\frac{aq^2+\lambda q^3}{1+bq^3}\fplus \cr
&\frac{\lambda q^4-abq^6}{1+bq^4}
\fplus \fdots \fplus 
\frac{aq^{s+1}+\lambda q^{2s+1}}{1+bq^{2s+1}}\fplus \frac{\lambda q^{2s+2}-abq^{3s+3}}
{1+bq^{2s+2}}
\fplus \fdots .
\label{G-cfrac4-2}
\end{align}
\end{ex}
Using Exercise \ref{ex3.2}, \S \ref{sec:RR}, it is easy to see that when $a\to 0$, this continued fraction reduces to \eqref{g-cfrac2}.  This is a special case of a continued fraction of Heine \cite{heine}. Surprisingly, it does not appear in Ramanujan's work, but then it is not as good-looking as  Ramanujan's own continued fraction expansions for $G(aq, b, \lambda q)/G(a, b, \lambda)$. 


\section{The role of transformation formulas}\label{sec:trans}
One of Ramanujan's own continued fractions for $G(aq, b, \lambda q)/G(a, b, \lambda)$ is \cite[Entry 6.2.1]{LN1}: 
\begin{equation}\label{G-cfrac1}
\frac{G(aq, b, \lambda q)}{G(a, b, \lambda)}
=\frac{1}{1}\fplus\frac{aq+\lambda q}{1}\fplus\frac{bq+\lambda q^2 }{1}\fplus\frac{aq^2+\lambda q^3}{1}\fplus\frac{bq^2+\lambda q^4}{1} \fdots ,
\end{equation}
where $G(a, b, \lambda)$ is defined in \eqref{G}.

Our strategy to derive \eqref{G-cfrac1} is to first transform the sums in \eqref{G-frac-sums} to obtain another ratio of similar looking sums. Now, using \eqref{div-1step}, the ratio of these transformed sums is expanded into Ramanujan's continued fraction \eqref{G-cfrac1}.

The transformation formula that we need was also known to Ramanujan 
(see {\cite[ch.~16, Entry 8]{berndt-notebooks-3}}): Let $|q|<1$ and $|a|<1$. Then
\begin{equation} \label{entry8}
\frac{\qrfac{a}{\infty}}{\qrfac{b}{\infty}}
\sum_{k=0}^{\infty} \frac{\qrfac{b/a}{k}\qrfac{c}{k}}{\qrfac{d}{k}\qrfac{q}{k}}a^k
=\sum_{k=0}^{\infty} \frac{\qrfac{b/a}{k}\qrfac{d/c}{k}}{\qrfac{b}{k}\qrfac{d}{k}\qrfac{q}{k}}
(ac)^k (-1)^k q^{\frac{k(k-1)}{2}}.
\end{equation}
In this transformation formula, we set $d=0$; $b\mapsto -bq$, $a\mapsto abq/\lambda$ and $c\mapsto -C$, to obtain:
\begin{align}
\sum_{k=0}^{\infty} \frac{\qrfac{-\lambda/a}{k}}{\qrfac{-bq}{k}\qrfac{q}{k}} &
\left(\frac{abC}{\lambda}\right)^k q^{\frac{k(k+1)}{2}}
=
\frac{\qrfac{abq/\lambda}{\infty}}{\qrfac{-bq}{\infty}}\cr
&\times
\sum_{k=0}^{\infty} \frac{\qrfac{-\lambda/a}{k}\qrfac{-C}{k}}{\qrfac{q}{k}}\left(\frac{abq}{\lambda}\right)^k.
 \label{entry8-d0}
\end{align}
Here, the RHS of \eqref{entry8} reduces to the LHS of \eqref{entry8-d0}.
Now notice that if $C=\lambda q/b$, we get $G(aq, b,\lambda q)$ on the LHS of \eqref{entry8-d0}. Similarly, set
$C=\lambda/b$ to obtain $G(a, b,\lambda )$. Taking ratios, we find that, for $| abq/\lambda |<1$:
\begin{equation}\label{G-frac-sums2}
\frac{G(aq, b, \lambda q)}{G(a, b, \lambda)}
=
 \frac{\displaystyle
\sum_{k=0}^{\infty}\frac{\qrfac{-\lambda/a}{k}\qrfac{-\lambda q/b}{k}}
{\qrfac{q}{k}}\left( \frac{abq}{\lambda}\right)^k}
{\displaystyle
\sum_{k=0}^{\infty}\frac{\qrfac{-\lambda/a}{k}\qrfac{-\lambda/b}{k}}
{\qrfac{q}{k}}\left( \frac{abq}{\lambda}\right)^k}.
\end{equation}
Expanding this ratio as a continued fraction, as usual using \eqref{div-1step} repeatedly, we can derive Ramanujan's continued fraction \eqref{G-cfrac1}. The calculations are very similar to the corresponding calculations for  \eqref{g-cfrac2} and Exercise \eqref{G-cfrac4-2}; if anything, they are simpler and more elegant.   Here is a compressed  proof.

Define $G_{1A}(s)$ and $G_{1B}(s)$, for $s=0,1, 2, 3, \dots$, as follows:
\begin{equation*}
G_{1A}(s):={\displaystyle\sum_{k=0}^{\infty} \frac{\qrfac{-\lambda q^s/a}{k} \qrfac{-\lambda q^s/b}{k}
}{\qrfac{q}{k}}\left(\frac{abq}{\lambda}\right)^k};
\end{equation*}
and,
\begin{equation*}
G_{1B}(s):={\displaystyle\sum_{k=0}^{\infty} \frac{\qrfac{-\lambda q^s/a}{k} \qrfac{-\lambda q^{s+1}/b}{k}
}{\qrfac{q}{k}}\left(\frac{abq}{\lambda}\right)^k}.
\end{equation*}
Then we find, for $s=0,1, 2, 3, \dots$
\begin{align*}
\frac{G_{1A}(s)}{G_{1B}(s)}&=
1+\frac{aq^{s+1}+\lambda q^{2s+1}}{\cfrac{G_{1B}(s)}{G_{1A}(s+1)}};\\
\frac{G_{1B}(s)}{G_{1A}(s+1)}&=
1+\frac{bq^{s+1}+\lambda q^{2s+2}}{\cfrac{G_{1A}(s+1)}{G_{1B}(s+1)}}.
\end{align*}
Now from \eqref{G-frac-sums2} and by iteration, we obtain:
\begin{align*}
&\frac{G(aq, b, \lambda q)}{G(a, b, \lambda)}
=
\frac{G_{1B}(0)}{G_{1A}(0)}=\frac{1}{G_{1A}(0)/G_{1B}(0)}
\cr
&=
\frac{1}{1}\fplus\frac{aq+\lambda q}{1}\fplus\frac{bq+\lambda q^2}
{1}\fplus\frac{aq^2+\lambda q^3}{1}\fplus \cr
&\frac{bq^2+\lambda q^4}{1}
\fplus \fdots \fplus 
\frac{aq^{s+1}+\lambda q^{2s+1}}{1}\fplus \frac{bq^{s+1}+\lambda q^{2s+2}}
{\cfrac{G_{1A}(s+1)}{G_{1B}(s+1)}}
.
\end{align*}
Ramanujan's continued fraction \eqref{G-cfrac1} follows by taking the limit as $s\to \infty$. Once again, we need to use the approach of \S\ref{sec:convergence} to justify the limiting process. 

Recall the definition \eqref{g} of $g(b,\lambda)$.  When $a\to 0$ in \eqref{G-cfrac1}, we obtain another continued fraction for $g(b, \lambda q)/g(b,\lambda)$ given by Ramanujan \cite[Entry 6.3.1 (ii)]{LN1}, that is different from \eqref{g-cfrac2}. 
\begin{equation}
\frac{g(b,\lambda q)}{g(b,\lambda)}
=
\frac{1}{1}\fplus\frac{\lambda q}{1}\fplus\frac{\lambda q^2 + bq}{1}\fplus\frac{\lambda q^3}{1}\fplus\fdots .\label{g-cfrac1}
\end{equation}

Now the exercises. 
\begin{ex}[Eisenstein (1844)]
Prove the following continued fraction due to Eisentein rediscovered by Ramanujan (see \cite[ch.~16, Entry 13]{berndt-notebooks-3} or \cite[Cor.~6.2.5]{LN1}):
\begin{equation}
\sum_{k=0}^{\infty} (-a)^k q^{\frac{k(k+1)}{2}}=
\frac{1}{1}\fplus\frac{aq}{1}\fplus\frac{a(q^2 - q)}{1}\fplus\frac{aq^3}{1}\fplus\frac{a(q^4-q^2)}{1}
\fplus\fdots .\label{entry13}
\end{equation}
\end{ex}
One can either use \eqref{div-1step} directly, or observe that the continued fraction in \eqref{g-cfrac1} reduces to \eqref{entry13} when $\lambda =a$ and $b=-a$. However, the rule is, that you can only use 
Ramanujan's entries from \cite[ch.~16]{berndt-notebooks-3} (or mentioned here) in order to get the LHS of \eqref{entry13}.
The case $a=1$ was known to Gauss in 1797, see \cite[p.~152]{LN1}.

The next exercise, requires the following transformation formula known to 
Ramanujan~\cite[ch.~16, Entry 6]{berndt-notebooks-3}: Let $|q|<1$, $|a|<1$, and $|c|<1$. Then
\begin{equation} \label{entry6}
\frac{\qrfac{a}{\infty}}{\qrfac{b}{\infty}}
\sum_{k=0}^{\infty} \frac{\qrfac{c}{k}\qrfac{b/a}{k}}{\qrfac{d}{k}\qrfac{q}{k}}a^k
=
\frac{\qrfac{c}{\infty}}{\qrfac{d}{\infty}}
\sum_{k=0}^{\infty} \frac{\qrfac{a}{k}\qrfac{d/c}{k}}{\qrfac{b}{k}\qrfac{q}{k}}c^k
.
\end{equation}
This is an easy consequence of \eqref{entry8}. 
\begin{ex} Use a special case of the transformation formula \eqref{entry6} to show that for $|\lambda/a|<1$:
\begin{align}
&\frac{G(aq, b, \lambda q)}{G(a, b, \lambda)}
=
 \frac{\displaystyle
\sum_{k=0}^{\infty}\frac{\qrfac{abq/\lambda}{k}}{\qrfac{q}{k}\qrfac{-aq^2}{k}}
\left(-\frac{\lambda}{a}\right)^k}
{(1+aq)\displaystyle\sum_{k=0}^{\infty}\frac{\qrfac{abq/\lambda}{k}}{\qrfac{q}{k}\qrfac{-aq}{k}}
\left(-\frac{\lambda}{a}\right)^k}
\label{G-frac5-sums}
\\
&=
\frac{1}{1+aq}\fplus\frac{\lambda q-abq^2}{1+aq^2}\fplus\frac{\lambda q^2+bq}
{1+aq^3}\fplus\frac{\lambda q^3-abq^5}{1+aq^4}\fplus \cr
&\frac{\lambda q^4+bq^2}{1+aq^5}
\fplus \fdots \fplus 
\frac{\lambda q^{2s+1}-abq^{3s+2}}{1+aq^{2s+2}}\fplus \frac{\lambda q^{2s+2}+bq^{s+1}}
{1+aq^{2s+3}}
\fplus \fdots .
\label{G-cfrac5}
\end{align}
\end{ex}
Note that we cannot take the limit as $a\to 0$ in the sums appearing on the RHS of \eqref{G-frac5-sums}. However, the continued fraction reduces to \eqref{g-cfrac2} as $a\to 0$. This continued fraction too does not appear in Ramanujan's work. However, it does appear to be related to both \eqref{G-cfrac4-2} and \eqref{G-cfrac1}. 

So far, we have seen two continued fractions 
for $g(b,\lambda q)/g(b,\lambda)$ (equations \eqref{g-cfrac1} and \eqref{g-cfrac2}), and three
for 
$G(aq, b, \lambda q)/G(a, b, \lambda)$ (equations \eqref{G-cfrac1}, \eqref{G-cfrac4-2} and \eqref{G-cfrac5}). However, that is not enough for an inventive genius like Ramanujan. He found more!

\section{A dose of insight into algebraical formulae}\label{sec:insight}

Ramanujan has yet another continued fraction for $g(b,\lambda q)/g(b,\lambda)$ in addition to \eqref{g-cfrac2} and \eqref{g-cfrac1}, see \cite[6.3.1 (iii)]{LN1}:
\begin{equation}
\frac{g(b,\lambda q)}{g(b,\lambda)}
=\frac{1}{1-b}\fplus\frac{b+\lambda q}{1-b}\fplus\frac{b+\lambda q^2}{1-b}\fplus\frac{b+\lambda q^3}{1-b}\fplus\fdots.\label{g-cfrac3}
\end{equation}
We will prove this formula under the additional requirement that $|b/(1-b)^2|<1/4$. 

To derive this formula, we use a transformation formula (a special case of the $a\to 0$ case of \eqref{entry8}) and obtain an equivalent ratio of sums. Once again, we use Euler's approach. 
But this time, we need a small trick to simplify the calculations a bit. This leads to \eqref{g-cfrac3}. This small dose of algebraical insight  is (I promise) the last trick required to derive such continued fractions. 

First, take limits as $a\to 0$ in \eqref{G-frac-sums2} to obtain:
\begin{equation}\label{g-frac-sums2}
\frac{g(b, \lambda q)}{g(b, \lambda)}
=
 \frac{\displaystyle
\sum_{k=0}^{\infty}\frac{\qrfac{-\lambda q/b}{k}}
{\qrfac{q}{k}}b^kq^{\frac{k^2+k}{2}}}
{\displaystyle
\sum_{k=0}^{\infty}\frac{\qrfac{-\lambda/b}{k}}
{\qrfac{q}{k}}b^kq^{\frac{k^2+k}{2}}}.
\end{equation}

We define $g_2(s)$, for $s=0,1, 2, 3, \dots$, as follows:
\begin{equation}\label{g3s}
g_2(s):={\displaystyle\sum_{k=0}^{\infty} \frac{\qrfac{-\lambda q^s/b}{k}}{\qrfac{q}{k}}b^k q^{\frac{k^2+k}{2}}}.
\end{equation}

Now, using \eqref{div-1step} we have, for $s=0,1, 2, 3,\dots$:
\begin{align*}\allowdisplaybreaks
\frac{g_2(s)}{g_2(s+1)}&=1+\frac{g_2(s)-g_2(s+1)}{g_2(s+1)}\cr
&=1+ \frac{1}{g_2(s+1)}
{\displaystyle\sum_{k=0}^{\infty} \frac{\qrfac{-\lambda q^s/b}{k}}{\qrfac{q}{k}}b^k q^{\frac{k^2+k}{2}}
\left[1-\frac{1+\lambda q^{s+k}/b}{1+\lambda q^{s}/b}\right]}
\cr
&=1+\frac{1}{g_2(s+1)}
{\sum_{k=0}^{\infty} \frac{\qrfac{-\lambda q^{s+1}/b}{k}}{\qrfac{q}{k}}b^k q^{\frac{k^2+k}{2}} 
\lambda q^{s+k+1}},
\end{align*}
 after canceling the first term in the difference of sums, and shifting the index to make the sum run from $0$ to $\infty$. 

 So far, the calculations are the same as before. But now its time for the algebraical insight mentioned earlier. Motivated by
the desire to get $g_2(s+2)$ in the sum on the RHS, and thus get a factor 
$\qrfac{-\lambda q^{s+2}/b}{k}$ in the sum, we use the following nice trick:
$$\lambda q^{s+k+1} = \lambda q^{s+k+1} +b-b = b(1+\lambda q^{s+k+1}/b)-b.$$
In addition, use
$$\qrfac{-\lambda q^{s+1}/b}{k}\cdot b(1+\lambda q^{s+k+1}/b)
=b(1+\lambda q^{s+1}/b) \qrfac{-\lambda q^{s+2}/b}{k}
$$
to obtain:
\begin{align*}
1+&\frac{1}{g_2(s+1)}
{\sum_{k=0}^{\infty} \frac{\qrfac{-\lambda q^{s+2}/b}{k}}{\qrfac{q}{k}}b^k q^{\frac{k^2+k}{2}} 
\left[b\left(1+\lambda q^{s+1}/b\right)\right]}\\
-&\frac{b}{g_2(s+1)}
{\sum_{k=0}^{\infty} \frac{\qrfac{-\lambda q^{s+1}/b}{k}}{\qrfac{q}{k}}b^k q^{\frac{k^2+k}{2}} 
}\\
=&
1+\frac{b+\lambda q^{s+1}}{g_2(s+1)}g_2(s+2)-\frac{b}{g_2(s+1)}\cdot {g_2(s+1)} \cr
=&
1-b+\frac{b+\lambda q^{s+1}}{g_2(s+1)/g_2(s+2)}.
\end{align*}
In this manner, we obtain the recurrence relation, for $s=0, 1, 2, 3, \dots$:
\begin{equation}\label{g3-recursion}
\frac{g_2(s)}{g_2(s+1)}=
1-b+\frac{b+\lambda q^{s+1}}{{\cfrac{g_2(s+1)}{g_2(s+2)}}}.
\end{equation}
This gives, on iteration,
\begin{equation}\label{g-cfrac2-finite}
\frac{g(b,\lambda q)}{g(b,\lambda)}=\frac{1}{g_2(0)/g_2(1)}=
\frac{1}{1-b}\fplus\frac{b+\lambda q}{1-b}\fplus \frac{b+\lambda q^2}{1-b}\fplus\fdots 
\fplus\frac{b+\lambda q^{s+1}}{\cfrac{g_2(s+1)}{g_2(s+2)}}.
\end{equation}
Note that \eqref{g-cfrac2-finite} implies modified convergence of the continued fraction. However, to be able to apply Worpitzky's Theorem to prove ordinary convergence, we need to consider the equivalent continued fraction
$$
\frac{1/(1-b)}{1}\fplus\frac{(b+\lambda q)/(1-b)^2}{1}\fplus \frac{\left( b+\lambda q^2\right)/(1-b)^2}{1}\fplus\fdots 
.
$$
Now assuming that  $|b/(1-b)^2|<1/4$, we can use the approach of 
\S\ref{sec:convergence} to prove ordinary convergence of Ramanujan's continued fraction \eqref{g-cfrac3}. 


Now some opportunities to develop your own insight into algebraical formulae. 
\begin{ex}
Define $G_2(s)$, for $s=0,1, 2, 3, \dots$, and $|\lambda/a|<1$ as follows:
\begin{equation}\label{G2s}
G_2(s):={\displaystyle\sum_{k=0}^{\infty} \frac{\qrfac{abq^s/\lambda}{k}}{\qrfac{q}{k}\qrfac{-aq^{s+1}}{k}}\left(-\frac{\lambda}{a}\right)^k }.
\end{equation}
Use \eqref{G-frac5-sums} to show that 
$$\frac{G(aq,b,\lambda q)}{G(a,b,\lambda )} =\frac{1}{1+aq}\fplus \frac{\lambda q-abq^2}
{(1+aq^2) \cfrac {G_2(1)}{G_2(2)}}.$$
Further, show that
$$
(1+aq^{s+1})\frac{G_{2}(s)}{G_{2}(s+1)}=
1+aq^{s+1}+bq^s+\frac{\lambda q^{s+1}-abq^{2s+2}}{(1+aq^{s+2})\cfrac{G_{2}(s+1)}{G_{2}(s+2)}},
$$
and derive Ramanujan's \cite[Entry 6.4.1]{LN1}  continued fraction:
\begin{align}
\frac{G(aq, b, \lambda q)}{G(a, b, \lambda)}
&=\cr
\frac{1}{1+aq}&\fplus\frac{\lambda q - abq^2}{1+aq^2+bq}\fplus\frac{\lambda q^2 -abq^4 }
{1+aq^3+bq^2}\fplus\frac{\lambda q^3-abq^6}{1+aq^4+bq^3}\fplus \fdots .\label{G-cfrac2}
\end{align}
\end{ex}

\begin{ex}[Hirschhorn \cite{mdh}, Bhargava and Adiga \cite{bhargava-adiga}] Use \eqref{G-frac-sums} to show that
\begin{equation}\label{G-cfrac3} 
\frac{G(aq, b, \lambda q)}{G(a, b, \lambda)}
=
\frac{1}{1}\fplus\frac{aq+\lambda q}{1-aq+bq}\fplus\frac{aq+\lambda q^2}
{1-aq+bq^2}\fplus\frac{aq+ \lambda q^3}{1-aq+bq^3}\fplus \fdots .
\end{equation} 
To show ordinary convergence, use the additional condition $|aq/(1-aq)^2|<1/4$.
\end{ex}
Note that when $b=0$, $a\mapsto b/q$, this reduces to Ramanujan's continued fraction \eqref{g-cfrac3}. See \cite[Theorem 6.4.1]{LN1} for a very similar proof of this continued fraction. 
This continued fraction appears in different forms in Hirschhorn \cite{mdh} and  Bhargava and Adiga \cite{bhargava-adiga}. We need to change a few parameters and appeal to \eqref{entry8-d0} in order to match one to the other.


So far, we have considered the general continued fractions of Ramanujan that appear in 
\cite[ch.~6]{LN1}. These consist of three continued fractions 
\eqref{g-cfrac1}, \eqref{g-cfrac2} and \eqref{g-cfrac3} for $g(b,\lambda q)/g(b,\lambda)$.  Further, Ramanujan noted two continued fractions, namely \eqref{G-cfrac1} and
\eqref{G-cfrac2}, for $G(aq, b, \lambda q)/G(a, b, \lambda)$. In addition, we mentioned three closely related continued fraction expansions: a continued fraction \eqref{G-cfrac3}  due to Hirschhorn \cite{mdh} (see also \cite{bhargava-adiga}); and two more, namely  \eqref{G-cfrac4-2} and \eqref{G-cfrac5}, that appear here for the first time. However, these last two cannot really be considered new, since they follow from a continued fraction of Heine \cite{heine}. 

Ramanujan wrote down many special cases of his general continued fractions. That is the subject of the next section. 

\section{Infinite products as continued fractions}\label{sec:infinite}

The most devilishly difficult of Ramanujan's formulas are often  particular cases of easier-to-derive general formulas.  
Consider, for instance,
Ramanujan's {\cite[Cor.~6.2.1]{LN1}} continued fraction formula, for $|q|<1$
\begin{align}
\frac{1}{1}\fplus\frac{q}{1}\fplus\frac{ q + q^2}{1}\fplus\frac{q^3}{1}\fplus
\frac{ q^2 + q^4}{1}\fplus\fdots 
&=
\frac{\pqrfac{q}{\infty}{q^2}}{(q^2; q^4)^{ 2}_{\infty}}\label{cor6.2.1-ln} \\
&=\frac{(1-q)(1-q^3)(1-q^5)\cdots }{(1-q^2)^2(1-q^6)^2(1-q^{10})^2\cdots }.\notag
\end{align}
Here the 
continued fraction is written as a ratio of infinite products, rather than as a ratio of sums. 
This happens when the sums themselves can be written as products. 

To prove \eqref{cor6.2.1-ln}, consider the $b=1$ and $\lambda =1 $ case of Ramanujan's continued fraction \eqref{g-cfrac1}. We find that the continued fraction equals $g(1, q)/g(1,1)$, where $g(b,\lambda)$ is defined in \eqref{g}. 
Observe that
$$\qrfac{q}{k}\qrfac{-q}{k}=\pqrfac{q^2}{k}{q^2}$$
and so
$$g(1,\lambda)=\sum_{k=0}^{\infty} \frac{q^{k^2}}{\pqrfac{q^2}{k}{q^2}}\lambda^k.$$
The sum $g(1, \lambda)$ can be written as an infinite product. For this we need the $q$-analog of the Binomial Theorem due to Rothe (1811): 
\begin{equation} \label{qbin}
\sum_{k=0}^{\infty} \frac{\qrfac{-b/a}{k}}{\qrfac{q}{k}}a^k
=
\frac{\qrfac{-b}{\infty}}{\qrfac{a}{\infty}},
\end{equation}
where,  $|q|<1$ and $|a|<1$.
The $q$-binomial theorem was known to Ramanujan. It is Entry 2 of \cite[Ch.~16]{berndt-notebooks-3}, and follows from \eqref{entry6} by setting $c=d$ and $b\mapsto -b$.

If we take the limit as $a\to 0$ in \eqref{qbin}, replace $q$ by $q^2$ and set $b\mapsto \lambda q$, we obtain
$$g(1,\lambda)=\pqrfac{-\lambda q}{\infty}{q^2}.$$
Thus, using some parts of the exercise in \S\ref{sec:RR},  we find that the continued fraction in \eqref{cor6.2.1-ln} equals
$$\frac{g(1,q)}{g(1,1)} = \pqrfac{-q^2}{\infty}{q^2} \times \frac{1}{\pqrfac{-q}{\infty}{q^2}}
=\frac{1}{\pqrfac{q^2}{\infty}{q^4}} \times
\frac{\pqrfac{q}{\infty}{q^2}}{\pqrfac{q^2}{\infty}{q^4}}
=\frac{\pqrfac{q}{\infty}{q^2}}{(q^2; q^4)^{ 2}_{\infty}},
$$ 
as required.

The final exercise outlines the derivation of another continued fraction of Ramanujan, found in Entry 11 of Chapter 16 of \cite{berndt-notebooks-3}. The original proof appears in \cite{abbw}. 

\begin{ex}
Show that, for $|q|<1$ and $|a|<1$
\begin{align*}
 \frac{\displaystyle
\sum_{k=0}^{\infty}\frac{\qrfac{b/a}{2k+1}}
{\qrfac{q}{2k+1}}a^{2k+1}}
{\displaystyle
\sum_{k=0}^{\infty}\frac{\qrfac{b/a}{2k}}
{\qrfac{q}{2k}}a^{2k}}
&=
 \frac{\displaystyle \frac{\qrfac{b}{\infty}}{\qrfac{a}{\infty}}-\frac{\qrfac{-b}{\infty}}{\qrfac{-a}{\infty}}}
{\displaystyle \frac{\qrfac{b}{\infty}}{\qrfac{a}{\infty}}+\frac{\qrfac{-b}{\infty}}{\qrfac{-a}{\infty}}}
\cr
&=
\frac{{\qrfac{-a}{\infty}}{\qrfac{b}{\infty}}-{\qrfac{a}{\infty}}{\qrfac{-b}{\infty}}}
{{\qrfac{-a}{\infty}}{\qrfac{b}{\infty}}+\qrfac{a}{\infty}{\qrfac{-b}{\infty}}}.
\end{align*}
Define, for $s=1, 2, 3, \dots$
$$C(s):= \sum_{k=0}^{\infty}\frac{\qrfac{bq^s/a}{2k}}
{\qrfac{q^2}{2k}}a^{2k}
\prod_{i=1}^{s-1} \frac{1-q^{2i+1}}{1-q^{2k+2i+1}}
.$$
With this definition, show that:
$$
\frac{\displaystyle
\sum_{k=0}^{\infty}\frac{\qrfac{b/a}{2k+1}}
{\qrfac{q}{2k+1}}a^{2k+1}}
{\displaystyle
\sum_{k=0}^{\infty}\frac{\qrfac{b/a}{2k}}
{\qrfac{q}{2k}}a^{2k}}
= \frac{a-b}{1-q}\fplus\frac{(a-bq)(aq-b)}{(1-q^3)\cfrac{C(1)}{C(2)}};
$$
and,
\begin{equation*}
(1-q^{2s+1})\frac{C(s)}{C(s+1)}=
1-q^{2s+1}+q^s\frac{(a-bq^{s+1})(aq^{s+1}-b)}{(1-q^{2s+3}){\cfrac{C(s+1)}{C(s+2)}}}.
\end{equation*}
Thus, for $|q|<1$ and $|a|<1$, derive the following continued fraction due to Ramanujan \cite[Entry 11, Ch. 16]{berndt-notebooks-3}:
\begin{align}
\frac{{\qrfac{-a}{\infty}}{\qrfac{b}{\infty}}-{\qrfac{a}{\infty}}{\qrfac{-b}{\infty}}}
{{\qrfac{-a}{\infty}}{\qrfac{b}{\infty}}+\qrfac{a}{\infty}{\qrfac{-b}{\infty}}}
&=\cr
\frac{a-b}{1-q}\fplus\frac{(a-bq)(aq-b)}{1-q^3}\fplus 
&\frac{q(a-bq^2)(aq^2-b)}{1-q^5}\fplus\fdots. \label{entry11}
\end{align}
\end{ex}

Ramanujan wrote down many special cases where continued fractions are written as products.  
These include the continued fraction that is the subject of \cite{absyz},  which is a special case of  \eqref{G-cfrac2}.
More such continued fractions can be found in \cite[ch.~6]{LN1} and 
\cite[ch.~16]{berndt-notebooks-3}. 
 More examples of Ramanujan-type continued fractions have been given by Gu and Prodinger \cite{gp}. 

This brings us to the end of our  tutorial.

\section{Conclusion}

Speaking of Ramanujan,
Hardy  \cite[p.~xxxv]{hardy} famously remarked:
\begin{quote}
{\small\sl
It was his insight into algebraical formulae, transformation of infinite series, and so forth, that was most amazing. On this side most certainly I have never met his equal, and I can compare him only with Euler and Jacobi.}
\end{quote}

Our study of Ramanujan's continued fractions illustrates Hardy's comments.  All the continued fractions in this article are derived  using Euler's approach. In \S \ref{sec:trans}, we felt the need for using transformations of infinite series, in addition to Euler's approach. And  in \S \ref{sec:insight}, we saw how some
algebraic insight leads to better-looking formulas. 

We have got a glimpse of Ramanujan's amazing gifts. 
I hope, dear reader,  that it is enough to make you feel like developing your own insight into algebraical formulae, transformation of infinite series, and so forth!

%
%
%
%
%

\subsection*{Acknowledgments} This tutorial was written in 2012 to celebrate the 125$^\text{th}$ year of Ramanujan's birth.

\end{document}